\input amstex
\documentstyle{amsppt}

\def\qed{\qquad \vrule height6pt width5pt depth0pt}
\def\z{{\bold Z}}

\def\q{{\bold Q}}

\def\oh{{\Cal O}}

\def\calt{{\Cal T}}
\def\k{{\Cal K}}
\def\l{{\Cal L}}

\def\calr{{\Cal R}}
\def\s{{\Cal S}}
\def\Chi{{\raise2pt\hbox{$\chi$}}}
\def\restrictedto#1{\big|\lower3pt\hbox{$\scriptstyle #1$}}
\def\cstar{$C^*$}

\def\inv{^{-1}}
\def\rank{\hbox{rank}\,}
\def\os{_{00}}
\input epsf         
\epsfverbosetrue    
\def\Xy{\leavevmode
 \hbox{\kern-.1em X\kern-.3em\lower.4ex\hbox{Y\kern-.15em}}}

\def\thmoo{Definition 1.1}   
\def\thma{Definition 1.2}   
\def\thmb{Lemma 1.3}        
\def\thmc{Lemma 1.4}        
\def\thmd{Proposition 1.5}  
\def\thme{Corollary 1.6}    
\def\thmi{Lemma 2.1}        
\def\thmj{Theorem 2.2}      

\def\bib#1{[#1]}
\def\blackadara{1}
\def\blackadarb{2}
\def\effroskaminker{3}
\def\eilersloring{4}
\def\kaplansky{5}
\def\kirchberg{6}
\def\kumjianpask{7}
\def\lin{8}
\def\loringa{9}
\def\neubuser{10}
\def\phillips{11}
\def\rordam{12}
\def\rosenbergschochet{13}
\def\spielberga{14} 
\def\spielbergb{15} 
\def\spielbergc{16} 
\def\spielbergd{17} 
\def\szymanski{18}  
\def\zhang{19}

\topmatter
\title Weak Semiprojectivity for Purely Infinite \cstar-algebras
\endtitle
\author Jack Spielberg \endauthor
\address Department of Mathematics and Statistics,
Arizona State University,
Tempe, AZ  85287-1804
\endaddress
\email jack.spielberg\@asu.edu\endemail
\abstract
We prove that a separable, nuclear, simple purely infinite,  
\cstar-algebra satisfying the universal coefficient theorem
 is weakly semiprojective if and only if 
its $K$-groups are direct sums of cyclic groups.
\endabstract
\endtopmatter
\document

\head Introduction \endhead

The first definition of semiprojectivity for \cstar-algebras was given 
by Effros and Kaminker in the context of noncommutative shape theory 
(\bib\effroskaminker).  A more restrictive definition was given by 
Blackadar in \bib\blackadara.  Loring introduced a third definition, 
which he termed {\it weak semiprojectivity\/}, in his investigations 
of stability problems for \cstar-algebras defined by generators and 
relations (\bib\loringa).  
Recently Neub\"{u}ser has introduced a slew of variants, 
the most important being what he called {\it asymptotic 
semiprojectivity\/} (\bib\neubuser).  Using the authors' initials to 
represent the above notions, the implications among them are:  B 
$\Rightarrow$ N $\Rightarrow$ EK, L.

All versions of semiprojectivity are of the following form:  
$*$-homorphisms into inductive limit \cstar-algebras can be lifted (in 
some sense) to a finite stage of the limit (the precise definitions 
may be found in section 1, and in the references). 
As a consequence, among the first 
(and easiest) examples for which semiprojectivity was established are 
the Cuntz-Krieger algebras.  This drew attention to the class of 
separable, nuclear, simple purely infinite 
\cstar-algebras, now commonly referred to as {\it Kirchberg\/} 
algebras (\bib\rordam).
 Kirchberg, and independently Phillips, have shown that in the 
 presence of the universal coefficient theorem,
$K$-theory is a complete invariant for Kirchberg algebras
 (\bib\kirchberg, \bib\phillips).  
Blackadar proved in \bib\blackadarb\ that for such algebras, finitely 
generated $K$-theory is necessary for semiprojectivity in the sense 
of \bib\effroskaminker.  He conjectured that for these algebras 
finitely generated $K$-theory is sufficient for semiprojectivity in 
the sense of \bib\blackadara, and proved this for the case of free 
$K_0$ and trivial $K_1$.  Szyma\'nski extended this to the case where 
$\rank K_1 \le \rank K_0$ (\bib\szymanski),
and in \bib\spielbergb\ semiprojectivity 
was proved whenever $K_1$ is free.
The conjecture remains open in the case that $K_1$ has torsion.
The methods used in all work 
on the conjecture rely upon explicit models for these algebras, 
constructed from 
directed graphs.  In another direction, Neub\"{u}ser used abstract 
methods to show that (for the algebras under consideration) finitely 
generated $K$-theory is equivalent to asymptotic semiprojectivity.

In this paper we study weak semiprojectivity for UCT-Kirchberg algebras. 
We prove that such an algebra is weakly semiprojective if and 
only if its $K$-groups are direct sums of cyclic groups.  The key 
difficulty lies in dealing with torsion in $K_1$, where we are 
forced to use tensor products of known semiprojectives.  
Semiprojectivity is badly behaved with respect to tensor products, 
and we rely on Neub\"{u}ser's result to get started.  Our 
contribution is thus in extending to the case where the $K$-theory is 
not finitely generated.  Another crucial technical aid is an 
alternate characterization of weak semiprojectivity, due to Eilers 
and Loring (\bib\eilersloring).

Our method of proof uses explicit models for the \cstar-algebras 
constructed from  a hybrid object which is partly a 
directed graph and partly a 2-graph (in the sense of \bib\kumjianpask).
The construction of this object, and the proof that it defines a 
UCT-Kirchberg algebra having the desired $K$-theory and given by 
suitable generators and relations, appears in \bib\spielbergc.

The outline of the paper is as follows.  In section 1 we prove 
the necessity in the main theorem.  This involves a kind of finite 
approximation property for abelian groups.  In section 2 we prove 
the main theorem.
During the final stages of writing an earlier draft of
 this paper we learned of Huaxin 
Lin's preprint \bib\lin, where the same theorem is proved by 
different means.
\smallskip
The figures in this paper were prepared with \Xy-pic.

\head 1.  Direct Sums of Cyclic Groups \endhead

The definition of weak semiprojectivity that follows is not Loring's 
original one, but was proved to be equivalent to it in 
\bib\eilersloring, Theorem 3.1.

\definition{\thmoo} The \cstar-algebra $A$ is called {\it weakly 
semiprojective\/} if given a \cstar-algebra $B$ with ideals $I_1\subseteq 
I_2\subseteq\cdots\subseteq I=\overline{\cup_k I_k}$, 
  a $*$-homomorphism $\pi:A\to B/I$, a finite 
set $M\subseteq A$, and $\epsilon>0$, there exists $n$ and a 
$*$-homomorphism $\phi:A\to B/I_n$ such that 
$$\Vert\pi(x)-\nu_n\circ\phi(x)\Vert<\epsilon\quad\text{ for } 
x\in M,$$ 
where $\nu_n:B/I_n\to B/I$ is the quotient map.

It is sometimes convenient to replace the increasing sequence of 
ideals by a directed family.
\enddefinition

We remark that if $M$ and $\epsilon$ are omitted, and it is required 
that $\pi=\nu_n\circ\phi$, then we recover Blackadar's definition of 
semiprojectivity.  Neub\"{u}ser's definition of asymptotic 
semiprojectivity can be obtained by omitting $M$ and $\epsilon$, and 
replacing $\phi$ by a point-norm continuous path $\phi_t$ such that 
for every $x\in A$, $\lim_t\Vert\pi(x)-\nu_n\circ\phi_t(x)\Vert=0$.

\definition{\thma} An abelian group $G$ has {\it Property $C$\/} (for {\it 
cyclic\/} --- see \thmd\ below)
if for every 
finite set $F\subseteq G$ there exists a finitely generated abelian 
group $K$, and homomorphisms $\alpha:G\to K$, $\beta:K\to G$ such 
that $\beta\circ\alpha(x)=x$ for all $x\in F$.
\enddefinition

\proclaim{\thmb} Let $A$ be a UCT-Kirchberg algebra.  If $A$ is weakly 
semiprojective, then $K_*(A)$ has Property $C$. \endproclaim

\demo{Proof} By \bib\kirchberg, $A=\overline{\cup A_n}$, 
$A_n\subseteq A_{n+1}$, where each $A_n$ is a UCT-Kirchberg algebra
 with finitely generated 
$K$-theory.  We modify the mapping telescope construction slightly 
(see, e.g., \bib\loringa).  Let
$$\eqalignno{B&=\left\{f\in C\bigl([0,1],A\bigr)\Bigm| f(t)\in A_n\hbox{\ for\ } 
t\ge {1\over n}\right\},\cr
J_n&=\left\{f\in B\Bigm| f\restrictedto{[0,1/n]}=0\right\},\cr
J&=\overline{\bigcup J_n}=\left\{f\in B\Bigm|f(0)=0\right\}.\cr
\noalign{\hbox{Then}}
B/J_n&\cong\left\{f\restrictedto{[0,1/n]}\Bigm|f\in B\right\},\cr
B/J&\cong A.\cr}$$
Now let $[x_1]$, $\ldots$, $[x_k]\in K_*(A)$.  Then by weak 
semiprojectivity there are $n$ and $\phi:A\to 
B/J_n$ such that
$$\left\Vert\phi(x_i)(0)-x_i\right\Vert<1,\quad 1\le i\le k.$$
Since $\phi(x_i)(1/n)$ is homotopic to $\phi(x_i)(0)$, we have that
$$[x_i]=\bigl[\phi(x_i)(1/n)\bigr],\quad1\le i\le k.$$
Let $\alpha:K_*(A)\to K_*(A_n)$ be given by
$$\alpha([x])=\bigl[\phi(x)(1/n)\bigr],$$
and let $\beta:K_*(A_n)\to K_*(A)$ be induced from the inclusion. 
Then $\beta\circ\alpha\bigl([x_i]\bigr)=[x_i]$ for $1\le i\le k$. \qed
\enddemo

\proclaim{\thmc}  Let $G$ be a countable abelian group with 
Property $C$.  Then $G/G_{tor}$ is free.\endproclaim
\demo{Proof} By \bib\kaplansky, Exercise 52, it suffices to show that every 
finite rank subgroup of $G/G_{tor}$ is free.  So let $H\subseteq 
G/G_{tor}$ be a subgroup of finite rank.  Put $\overline 
H=\pi\inv(H)$, where $\pi$ is the quotient map of $G$ onto 
$G/G_{tor}$.  Let $e_1$, $\ldots$, $e_r$ be a basis for $H$.  Then we 
may write $\z^r\subseteq H\subseteq \q^r$ (relative to this basis).  
Let $\overline{e_1}$, $\ldots$, $\overline{e_r}\in\overline H$ with 
$\pi(\overline{e_i})=e_i$, $1\le i\le r$.  Let $K$, $\alpha:G\to K$ and 
$\beta:K\to G$ be as in Property $C$, with 
$\beta\circ\alpha(\overline{e_i})=\overline{e_i}$ for $1\le i\le r$.  

We claim that $\ker\bigl(\alpha\restrictedto{\overline{H}}\bigr)\subseteq 
G_{tor}$.  To see this, let $y\in
\ker\bigl(\alpha\restrictedto{\overline{H}}\bigr)$.  Choose $N\in\z$ 
such that $N\pi(y)\in\z^r$.  We may write
$$\pi(Ny)=\sum_{i=1}^rc_ie_i,\quad c_i\in\z.$$
Then
$$z=Ny-\sum_{i=1}^rc_i\overline{e_i}\in\ker\pi=G_{tor}.$$
Thus
$$0=N\beta\circ\alpha(y)=\beta\circ\alpha(Ny)=\beta\circ\alpha (z) + 
\sum_{i=1}^rc_i\overline{e_i}.$$
But since $\beta\circ\alpha(G_{tor})\subseteq G_{tor}$, we may 
apply $\pi$ to the last equation to get
$$0=\sum_{i=1}^rc_ie_i.$$
It follows that $c_i=0$ for all $i$, so that $Ny=z\in G_{tor}$.  
Hence $y\in G_{tor}$.

Next we claim that $\ker\bigl(\pi_K\circ\alpha\restrictedto{\overline 
H}\bigr)=G_{tor}$, where $\pi_K$ is the quotient map of $K$ onto 
$K/K_{tor}$.
To see this, first note that the containment 
$\supseteq$ is obvious.  For the other containment, let 
$y\in\ker\bigl( \pi_K\circ\alpha\restrictedto{\overline H}\bigr)$. 
Then $\alpha(y)\in K_{tor}$, so $N\alpha(y)=0$ for some 
$N\in\z\setminus\{0\}$.  Then $Ny\in\ker\bigl( 
\alpha\restrictedto{\overline H}\bigr)$, so $Ny\in G_{tor}$ by the 
previous claim.  Hence $y\in G_{tor}$.

Finally, it follows from the last claim that $\pi_K\circ\alpha 
\restrictedto{\overline H}$ induces an injection $H\to K/K_{tor}$, 
which implies that $H$ is free. \qed
\enddemo

\proclaim{\thmd} Let $G$ be a countable abelian group.  Then $G$ has
Property $C$ if and only if
$G$ is a direct sum of cyclic groups. \endproclaim

\demo{Proof} It is clear that a direct sum of cyclic groups has Property C.  
Conversely,
by \thmc, $G\cong G_{tor}\oplus G/G_{tor}$, where $G/G_{tor}$ is 
free, and hence a direct sum of (infinite) cyclic groups.  Since 
$G_{tor}=\oplus_p G_p$, where $G_p$ is the $p$-primary component of 
$G_{tor}$, it suffices to prove that $G_p$ is a direct sum of cyclic 
groups.  By \bib\kaplansky, Theorem 11, it suffices to prove that 
$G_p$ contains no element of infinite height.  To see this, let $x\in 
G_p\setminus\{0\}$.  Choose $K$, $\alpha:G\to K$, and $\beta:K\to  G$
as in Property $C$ so that $\beta\circ\alpha(x)=x$.  We have 
$\alpha(x)\in K_p$, the $p$-primary component of $K$.  Let $n$ be the 
maximal height of elements of $K_p$.  Now if $x=p^jy$ in $G$, then
$$\eqalign{x&=\beta\circ\alpha(x)\cr
&=\beta\circ\alpha(p^jy)\cr
&=\beta\bigl(p^j\alpha(y)\bigr)\cr
&=0,\quad\hbox{if\ }j> n.\cr}$$
Therefore $j\le n$, and so $x$ is of finite height. \qed
\enddemo

\proclaim{\thme}  Let $A$ be a UCT-Kirchberg algebra.  If $A$ is weakly 
semiprojective, then $K_*(A)$ is a direct sum of cyclic groups.
\qed\endproclaim

\head 2.  The Main Theorem \endhead

We now wish to prove the converse of \thme, establishing weak 
semiprojectivity for any
UCT-Kirchberg algebra whose $K$-theory is a direct 
sum of cyclic groups.  To do this we will use the models for 
UCT-Kirchberg algebras constructed in \bib\spielbergc.  For each 
$k\ge2$ let $H_k$ be the directed graph shown in figure 1.

\bigskip
\epsfysize=1.3 true in                   
\centerline{\epsfbox{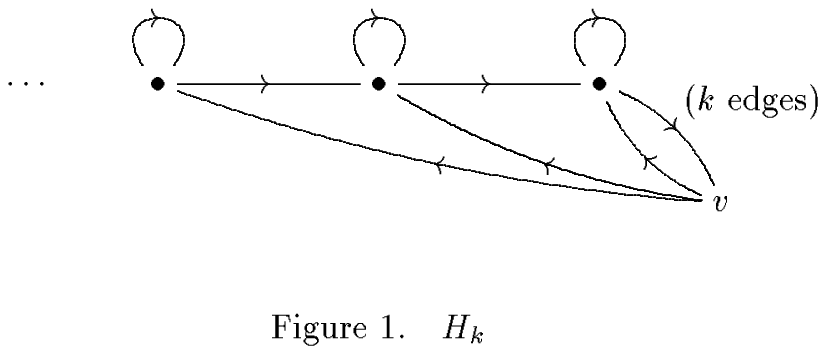}}
\bigskip

\noindent
One easily checks that $K_*\oh(E_k)=(\z/(k),0)$ (see e.g.
\bib\spielbergd). We let $H_\infty$ denote the usual directed graph 
of the Cuntz algebra $\oh_\infty$:  one vertex with denumerably many 
loops.  Finally we let $\overline{H}_\infty$ denote the graph shown in 
figure 2.

\bigskip
\epsfysize=1.3 true in                   
\centerline{\epsfbox{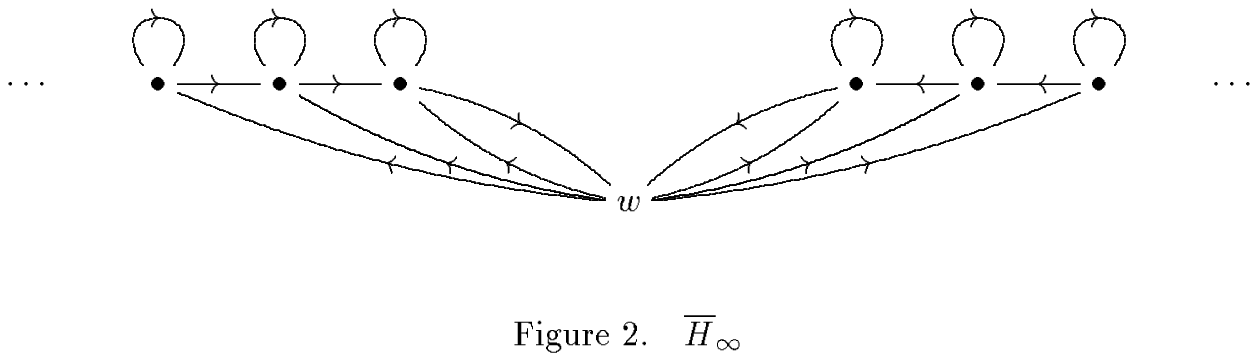}}
\bigskip

\noindent
Again one can easily check that 
$K_*\oh(\overline{H}_\infty)=(0,\z)$.  We remark that the graphs 
$H_k$, $H_\infty$ and $\overline{H}_\infty$ have a distinguished
vertex emitting infinitely many edges, as required for the
construction in \bib\spielbergc.

We now construct a UCT-Kirchberg algebra having as $K$-theory a 
prescribed direct sum of cyclic groups.  Let 
$G^i=(G^i_0, G^i_1)$ for $i=0$, $1$, $\ldots$, where for each $i$, 
one of $G^i_0$, $G^i_1$ is a cyclic group
and the other is the zero 
group.  By the K\"{u}nneth formula (\bib\rosenbergschochet),
$G^i$ is equal either to
$K_*\bigl(\oh(H_{k_i})\otimes\oh(H_\infty)\bigr)$
or to $K_*\bigl(\oh(H_{k_i})\otimes\oh(\overline H_\infty)\bigr)$
for some $k_i\in\{2,3,\ldots,\infty\}$.
Let $E_i=H_{k_i}$ and $F_i=H_\infty$ or $\overline{H}_\infty$
so that $G^i=K_*\bigl(\oh(E_i)\otimes\oh(F_i)\bigr)$.  As in 
\bib\spielbergc, we let $\Omega$ denote the hybrid object constructed 
from the product 2-graphs $E_i\times F_i$ and the connecting 
1-graphs $D_i$.  By Theorem 4.8 of \bib\spielbergc, $C^*(\Omega)$ 
is a UCT-Kirchberg algebra with $K$-theory equal to $\oplus_i G^i$.  
Let $A=C^*(\Omega)$.

We briefly recall the definition of the \cstar-algebra $C^*(\Omega)$
from \bib\spielbergc.  First let us recall the definitions of the
\cstar-algebras of a directed graph in a form convenient for this
purpose.
A {\it directed graph\/} $E$ consists of two sets, $E^0$ (the {\it 
 vertices\/}) and $E^1$ (the {\it edges\/}), together with two maps 
 $o$, $t:E^1\to E^0$ ({\it origin\/} and {\it terminus\/}).
 We let $\oh(E)$ denote the \cstar-algebra of $E$.  It 
 is the universal \cstar-algebra defined by generators $\bigl\{P_a 
 \bigm| a\in E^0\bigr\}$ and $\bigl\{S_e \bigm| e\in E^1\bigr\}$ with 
 the {\it Cuntz-Krieger relations\/}:
\roster
 \item"$\bullet$" $\bigl\{P_a \bigm| a\in E^0\bigr\}$ are pairwise 
 orthogonal projections.
 \item"$\bullet$" $S_e^*S_e=P_{t(e)}$, for $e\in E^1$.
 \item"$\bullet$" $o(e)=o(f) \Longrightarrow S_eS_e^* + S_fS_f^* \le 
 P_{o(e)}$, for $e$, $f\in E^1$ with $e\not=f$.
 \item"$\bullet$" $0 < \#\ E^1(a)<\infty \Longrightarrow P_a= 
 \sum\bigl\{S_eS_e^*\bigm| o(e)=a\bigr\}$, for $a\in E^0$,
\endroster
 \par\noindent
 where in the fourth relation we use the notation $E^1(a)$ to denote
 the set of edges with origin $a$.
 (These are a variant of the relations given in \bib\spielberga, 
 Theorem 2.21.)

 The relationship between the \cstar-algebras of a graph and a 
 subgraph are crucial to our methods.  We refer to 
 \bib\spielberga.  The results are as follows.  
 Let $E$ be a graph and let $F$ be a subgraph of $E$.  We let 
 $S=S(F)$ be 
 the set of vertices in $F^0$ that do not emit more edges in $E$ than 
 in $F$.  We let $\calt\oh(F,S)$ denote the relative Toeplitz 
 Cuntz-Krieger algebra of $F$ in $E$.  It is the universal 
 \cstar-algebra defined by generators $\bigl\{P_a 
 \bigm| a\in F^0\bigr\}$ and $\bigl\{S_e \bigm| e\in F^1\bigr\}$ with 
 the relations (as above) for $\oh(F)$, modified by requiring the 
 fourth relation only if $a\in S$. Then $\calt\oh(F,S)$ is the 
 \cstar-subalgebra of $\oh(E)$ generated by the projections and partial 
 isometries associated to the vertices and edges of $F$ 
 (\bib\spielberga, Theorem 2.35)
 
The hybrid object $\Omega$ is constructed from a directed graph $D$
(see figure 3), and the sequence of product 2-graphs $E_i\times F_i$
(\bib\kumjianpask).

\bigskip
\epsfysize=2.5in                   
\centerline{\epsfbox{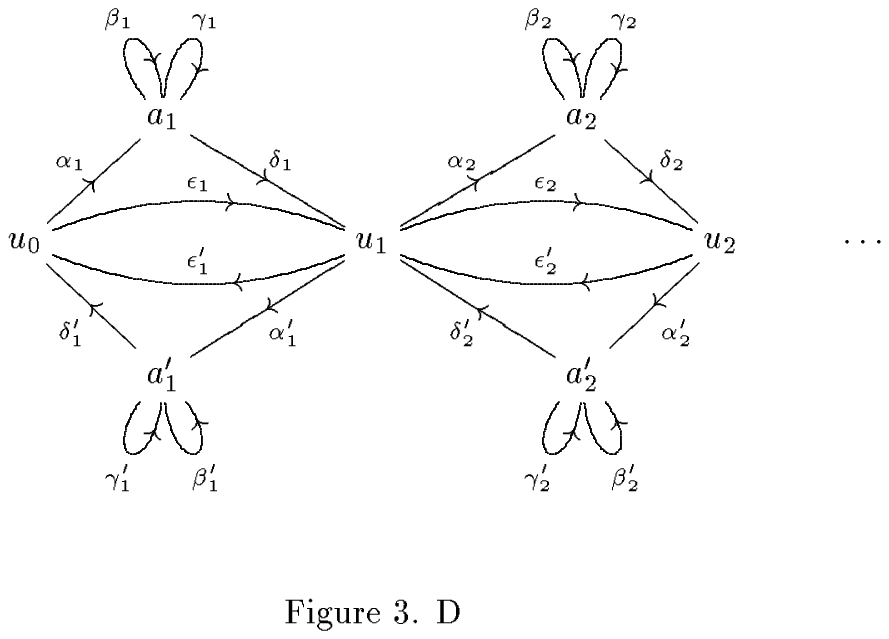}}
\bigskip

We let $v_i$, respectively $w_i$, denote the distinguished vertex in
$E_i$, respectively $F_i$, emitting infinitely many edges, and
we form $\Omega$ by attaching
$E_i\times F_i$ to $D$ by identifying $u_i$ with $(v_i,w_i)$.
By a {\it vertex\/} of $\Omega$ we mean an element of 
$$\cup_i (E_i^0\times F_i^0) \cup  
D^0,$$ 
where we identify $u_i$ and $(v_i,w_i)$.  By an {\it edge\/} we mean 
an element of $$\bigl(\cup_i(E_i^1\times F_i^0) \cup (E_i^0 
\times F_i^1)\bigr) \cup D^1.$$
The \cstar-algebra of $\Omega$ is defined by generators and relations 
as follows.
We let $\s$ 
denote the  set of symbols 
$$\bigl\{P_x\bigm|\ x \text{ is a vertex}\bigr\} \cup
\bigl\{S_y\bigm|\ y \text{ is an edge}\bigr\}.$$
We let $\calr$ denote the following set of relations on $\s$:
\smallskip\noindent
\item{(i)} $P_x$ is a projection for every vertex $x$,  $S_y$ is a 
partial isometry for every edge $y$.
\item{(ii)} For every $a\in E^0_i$, the projections for $\{a\}\times 
F_i^0$ and the partial isometries for $\{a\}\times F_i^1$ satisfy the 
Cuntz-Krieger relations corresponding to the graph $F_i$ (see 
the discussion at the end of the section 1).
\item{(ii')} For every $b\in F^0_i$, the projections for 
$E_i^0\times \{b\}$ and the partial isometries for 
$E_i^1\times \{b\}$ satisfy the 
Cuntz-Krieger relations corresponding to the graph $E_i$.
\item{(iii)} The projections for $D^0$ and 
the partial isometries for $D^1$ satisfy the Toeplitz-Cuntz-Kriger
relations corresponding to the graph $D$ and the vertices
$\{a_0,a_1\}$.
\item{(iv)} If $\mu$ and $\nu$ are edges of types $D$ 
and $E_i\times F_i$, respectively, then $S_\mu^* S_\nu=0$.
\item{(v)} For all $e\in E_i^1$ and $f\in F_i^1$ we have
$$\align
S_{(o(e),f)}\,S_{(e,t(f))} &= 
S_{(e,o(f))}\,S_{(t(e),f)} \\
S_{(t(e),f)}\,S_{(e,t(f))}^* &= 
S_{(e,o(f))}^*\, S_{(o(e),f)}.
\endalign$$
Then $A=C^*(\Omega)=C^*\langle\s,\calr\rangle$ 
is the universal \cstar-algebra given by these 
generators and relations.

\proclaim\thmi The $C^*$-algebra $A=C^*(\Omega)$ is weakly semiprojective.
\endproclaim

\demo{Proof}  Let $\Omega_{(n)}$ be the subobject of $\Omega$ 
consisting of $D_0$, $D_1$, $\ldots$, $D_n$, $E_0\times F_0$, 
$\ldots$, $E_n\times F_n$.  (We use parentheses in order to avoid 
confusion with the notation of \bib\spielbergc.)
Theorem 4.8 of \bib\spielbergc\ applies to 
$\Omega_{(n)}$, so that $A_n=C^*(\Omega_{(n)})$ is a UCT-Kirchberg 
algebra with finitely generated $K$-theory.  Note that the generators 
and relations defining $C^*(\Omega_{(n)})$ 
are the same in $C^*(\Omega_{(n)})$ as in 
$C^*(\Omega)$ (this is essentially because of the infinite valence of 
the vertex $u_{n+1}$).  Thus $A_n$ is a \cstar-subalgebra of $A$, 
and $A=\overline{\cup_n A_n}$.
By Satz 6.12 of \bib\neubuser, $A_n$ is uniformly asymptotically 
semiprojective.  It follows from Theorem 3.1 of \bib\eilersloring\ 
that $A_n$ is weakly semiprojective.

Let $B$ be a \cstar-algebra with ideals $I_1\subseteq 
I_2\subseteq\cdots\subseteq I=\overline{\cup_k I_k}$, and let 
$\pi:A\to B/I$ be a $*$-homomorphism.  Let $M\subseteq A$ be a finite 
set, and let $\epsilon>0$.  Choose $n$, and a finite set $M'\subseteq 
A_{n-1}$, such that $d(x,M')<\epsilon/2$ for all $x\in M$.  Since 
$A_n$ is weakly semiprojective there is $k$, and a $*$-homorphism 
$\phi_0:A_n\to B/I_k$, such that 
$$\Vert\pi(x')-\nu_k\circ\phi_0(x')\Vert<\epsilon/2\quad\text{ for } 
x'\in M',$$ 
where $\nu_k:B/I_k\to B/I$ is the quotient map.  We will construct 
a $*$-homomorphism $\phi:A\to B/I_k$ extending 
$\phi_0\restrictedto{A_{n-1}}$.  Then it will follow that
$$\Vert\pi(x)-\nu_k\circ\phi(x)\Vert<\epsilon\quad\text{ for }x\in M,$$
concluding the proof.

Let $p=P_{u_{n+1}}$ and $q=P_{u_n}$, the projections in $A$ 
corresponding to the vertices $u_{n+1}$ and $u_n$.  The hereditary 
subalgebra $pAp$ of $A$ contains a hereditary subalgebra, $C$, 
isomorphic to $A$.  (This follows easily from the pure infiniteness 
of $A$.  See, e.g., the proof of Theorem 3.12 in \bib\spielbergb.)  
Let $\psi_1:A\to C$ be a $*$-isomorphism.  Since the inclusion of $C$ 
into $A$ induces the identity in $K$-theory, it follows that $\psi_{1*}$ 
is an automorphism of $K_*(A)$.  It follows from Theorem 4.2.1 of 
\bib\phillips\ that there is a $*$-automorphism, $\alpha$, of $A$ 
with $\alpha_*=\psi_{1*}$.  Let $\psi_2=\psi_1\circ\alpha\inv$.  
Then $\psi_2:A\to C$ is a $*$-isomorphism, and $\psi_{2*}$ is the 
identity in $K$-theory.  Let $x\in A$ be a partial isometry with 
$x^*x=q$ and $xx^*=\psi_2(q)$.  Increasing $k$ if necessary we may 
find a partial isometry $z\in B/I_k$ with $z^*z=\phi_0(q)$ and 
$zz^*=\phi_0\circ\psi_2(q)$.  We define $\phi:A\to B/I_k$ by defining 
it on the generators $\s$ of $A$ (Definition 3.3 of \bib\spielbergc):
$$\phi(s_y)=\cases
\phi_0(s_y),& y\in\Omega_{n-1}\\
\phi_0\circ\psi_2(s_y),& y\not\in\Omega_{n-1}\text{ and }o(y),\;t(y) 
\not= u_n\\
\bigl(\phi_0\circ\psi_2(s_y)\bigr)z^*,& y\not\in\Omega_{n-1}\text{ and }t(y) 
= u_n,\ o(y)\not=u_n\\
z\bigl(\phi_0\circ\psi_2(s_y)\bigr),& y\not\in\Omega_{n-1}\text{ and }o(y) 
= u_n,\ t(y)\not=u_n\\
z\bigl(\phi_0\circ\psi_2(s_y)\bigr)z^*,& y\not\in\Omega_{n-1}\text{ and }o(y) 
= t(y) = u_n.\\\endcases$$
It is easy to see that the 
elements $\phi(s_y)$ satisfy the relations $\calr$ of 
\bib\spielbergc, Theorem 3.3, and 
hence $\phi$ defines a $*$-homomorphism. \qed
\enddemo

\proclaim{\thmj}  Let $A$ be a UCT-Kirchberg algebra.
  Suppose that $K_*(A)$ is a 
direct sum of cyclic groups.  
Then $A$ is weakly semiprojective. \endproclaim

\demo{Proof} As in the proof of Theorem 3.12 of \bib\spielbergb, it suffices 
to prove that if $A$ is unital and $\k\otimes A$ is weakly 
semiprojective, then $A$ is weakly semiprojective.
(We are relying on the classification theory of \bib\kirchberg\ and 
\bib\phillips, as well as the theorem of \bib\zhang\ that nonunital 
separable, simple, purely infinite \cstar-algebras are stable.)
Let $u_1$, $u_2$, $\ldots\in A$ with $u_i^*u_j=\delta_{ij}$.  Put
$$A_0=\overline{\hbox{span}}\,\bigl\{u_iAu_j^*\bigr\}.$$  Then 
$A_0$ is isomorphic to $\k\otimes A$.

Let $A=B/I$, where $I$ is the closure of a directed family 
of ideals, $\l$, of $B$.  (Since $A$ is simple we may dispense with 
the homomorphism $\pi$ of \thmoo.)
We let $\pi:B\to B/I$, $\pi_J:B\to B/J$ for 
$J\in\l$, denote the quotient maps.  Put $B_0=\pi\inv(A_0)$.  We will 
use \bib\eilersloring, Theorem 3.1.  So let $F\subseteq A$ be a 
finite set, and let $\epsilon>0$.  We may assume that $\epsilon<1$ and 
that $1\in F$.  Choose $\gamma<1$ such that $3\gamma\Vert 
x\Vert<\epsilon$ for all $x\in F$.
Since $A_0$ is weakly semiprojective by hypothesis, 
there is a $*$-homomorphism $\psi\os:A_0\to B/J$ such that
$$\Vert\pi\circ\psi\os(x)-x\Vert<\gamma\quad\hbox{for\ } x\in u_1 Fu_1^*.$$
In particular, we have
$$\Vert\pi\circ\psi\os(u_1u_1^*)-u_1u_1^*\Vert<\gamma.$$
Choose $v\in B$ with $\pi(v)=u_1$.  By increasing $J$, if necessary, we 
may assume that $\pi_J(v)$ is an isometry.  Since
$$\Vert\pi(vv^*)-\pi\circ\psi\os(u_1u_1^*)\Vert<\gamma,$$
we may assume, again by increasing $J$ if necessary, that
$$\Vert\pi_J(vv^*)-\psi\os(u_1u_1^*)\Vert<\gamma.$$
Then there exist $s$, $t\in B$ such that
$$\eqalign{\pi_J(s^*s)&=\pi_J(vv^*)\cr
\pi_J(ss^*)&=\psi\os(u_1u_1^*)\cr
\pi_J(t^*t)&=1-\pi_J(vv^*)\cr
\pi_J(tt^*)&=1-\psi\os(u_1u_1^*)\cr
\Vert\pi_J(s)-\pi_J(vv^*)\Vert&<\gamma\cr
\Vert\pi_J(t)-\bigl(1-\pi_J(vv^*)\bigr)\Vert&<\gamma.\cr}$$
Let $z=s+t$.  Then $\pi_J(z)$ is unitary, and 
$\Vert\pi_J(z)-1\Vert<\gamma$.  Define $\psi_0:A_0\to B/J$ by
$$\psi_0(x)=\pi_J(z)^*\psi\os(x)\pi_J(z).$$
Then $\psi_0$ is a $*$-homomorphism, and
$$\leqalignno{\Vert\pi\circ\psi_0(x)-x\Vert 
&=\Vert\pi(z)^*\pi\circ\psi\os(x)\pi(z)-x\Vert\cr
&\le2\Vert\pi(z)-1\Vert\,\Vert x\Vert + \Vert\pi\circ\psi\os(x)-x\Vert\cr
&\le2\gamma\Vert x\Vert+\gamma\cr
&<\epsilon,\quad\hbox{for\ } x\in 
u_1Fu_1^*,&(*)\cr
\psi_0(u_1u_1^*)&=\pi_J(z)^*\psi\os(u_1u_1^*)\pi_J(z)\cr
&=\pi_J(s)^*\psi\os(u_1u_1^*)\pi_J(s)\cr
&=\pi_J(vv^*).&(**)\cr}$$
Now define $\psi:A\to B/J$ by
$$\psi(x)=\pi_J(v)^*\psi_0(u_1xu_1^*)\pi_J(v).$$
By $(**)$ we have that $\psi$ is a $*$-homomorphism.  For $x\in F$, 
we have
$$\eqalign{\Vert\pi\circ\psi(x)-x\Vert&=\Vert 
u_1^*\pi\circ\psi_0(u_1xu_1^*) 
u_1 - x\Vert\cr
&=\Vert u_1^*\bigl(\pi\circ\psi_0(u_1xu_1^*) - u_1xu_1^*\bigr)u_1\Vert\cr
&\le\Vert\pi\circ\psi_0(u_1xu_1^*) - u_1xu_1^*\Vert\cr
&<\epsilon,\quad\hbox{by\ }(*).\qed\cr}$$
\enddemo

\head References \endhead

\roster
\item"\bib\blackadara" B. Blackadar, Shape theory for \cstar-algebras, {\it Math. 
Scand.\/} {\bf 56} (1985), 249-275.
\smallskip
\item"\bib\blackadarb"  B. Blackadar, Semiprojectivity in simple  
\cstar-algebras,  {\it Operator algebras and applications\/},  1--17, 
Adv. Stud. Pure Math., 38, Math. Soc. Japan, Tokyo, 2004.
\smallskip
\item"\bib\effroskaminker " E.G. Effros and J. Kaminker, Homotopy continuity and 
shape theory for \cstar-algebras, in {\it Geometric Methods in Operator 
Algebras\/}, eds. Araki and Effros, Pitman Res. Notes Math. {\bf 123},
Longman, Harlow, 1986.
\smallskip
\item"\bib\eilersloring " S. Eilers and T. Loring, Computing contingencies 
for stable relations, {\it Internat. J. Math.\/}
{\bf 10} (1999), no. 3, 301--326. 
\smallskip
\item"\bib \kaplansky " I. Kaplansky, {\it Infinite Abelian Groups\/},
The University of Michigan Press, Ann Arbor, 1969.
\smallskip
\item"\bib\kirchberg " E. Kirchberg, The classification of purely 
infinite \cstar-algebras using Kasparov's theory, {\it Fields 
Institute Communications\/}, 2000.
\smallskip
\item"\bib\kumjianpask " A. Kumjian and D. Pask, HIgher rank graph 
\cstar-algebras, {\it New York J. Math.\/} {\bf 6} (2000), 1-20.
\smallskip
\item"\bib\lin " H. Lin, Weak semiprojectivity in purely infinite 
simple \cstar-algebras, pre\-print (2002).
\smallskip
\item"\bib\loringa " T. Loring, {\it Lifting Solutions to Perturbing 
Problems in \cstar-algebras\/}, Fields Institute Monographs, vol. 8, 
Amer. Math. Soc., Providence, 1997.
\smallskip
\item"\bib\neubuser " B. Neub\"{u}ser, Semiprojektivit\"{a}t und realisierunen 
von rein unendlichen \cstar-alge\-bren, preprint, M\"{u}nster, 2000.
\smallskip
\item"\bib\phillips " N.C. Phillips, A classification theorem for 
nuclear purely infinite simple \cstar-algebras, {\it Document math.\/} 
(2000), no. 5,  49-114.
\smallskip
\item"\bib\rordam" M. R\o rdam, Classification of Nuclear
\cstar-Algebras, in {\it Encyclopaedia of Mathematical Sciences\/}
{\bf 126} Springer-Verlag (2002), Berlin Heidelberg New York.
\smallskip
\item"\bib\rosenbergschochet " J. Rosenberg and C. Schochet, The 
K\"{u}nneth theorem and the universal coefficient theorem for 
Kasparov's generalized $K$-functor, 
{\it Duke Math. J.\/} {\bf 55} (1987) no. 2, 431-474.
\smallskip
\item"\bib\spielberga " J. Spielberg, A functorial approach to the 
\cstar-algebras of a graph, {\it Internat. J. Math.\/} {\bf 13} 
(2002), no.3, 245-277.
\smallskip
\item"\bib\spielbergb " J. Spielberg, Semiprojectivity for certain 
purely infinite \cstar-algebras, pre\-print, Front for the
Mathematics ArXiv, math.OA/0102229.
\smallskip
\item"\bib\spielbergc" Graph-based models for Kirchberg algebras,
preprint (2005).
\smallskip
\item"\bib\spielbergd" Non-cyclotomic presentations of modules and 
prime-order auto\-mor\-phisms of Kirchberg algebras, preprint (2005).
\smallskip
\item"\bib\szymanski " W. Szyma\'nski, On semiprojectivity of \cstar-algebras of 
directed graphs, {\it Proc. AMS\/} {\bf 130} (2002), 1391-1399.
\smallskip
\item"\bib\zhang " S. Zhang, Certain \cstar-algebras with real rank zero and 
their corona and multiplier algebras, {\it Pacific J. Math.\/} {\bf 
155} (1992), 169-197.
\endroster

\enddocument
\bye